\documentclass[final,1p,times]{elsarticle}
\usepackage{graphicx}
\usepackage{amssymb}
\usepackage{amsthm}

\usepackage{setspace}
\usepackage[utf8]{inputenc}
\usepackage{algorithm}
\usepackage[noend]{algpseudocode}

\usepackage{enumitem} 
\usepackage{tablefootnote}
\usepackage{graphicx}
\usepackage{subfigure}
\usepackage{booktabs}
\usepackage{amsfonts}
\usepackage{amsthm}
\usepackage{amsmath,bm}
\usepackage{amssymb}
\usepackage{mathtools}
\usepackage{multirow}
\usepackage{bbm}
\usepackage{hyperref}

\usepackage{array}

\usepackage[colorinlistoftodos]{todonotes}
\usepackage{comment}

\usepackage{mathrsfs}

\usepackage{epstopdf}    
\usepackage{enumitem}  

\usepackage{cuted}%%\stripsep-3pt

\usepackage{utfsym}
\usepackage{diagbox}

\usepackage{float}
\usepackage{threeparttable}
\usepackage{url} % to cite web url in the bibtex

\DeclareMathOperator{\EX}{\mathbb{E}}% expected value

\usepackage{geometry}
\begin{document}

\begin{frontmatter}


% \ead{xiaoming.li@concordia.ca, hubert.normandin-taillon@mail.concordia.ca, chun.wang@concordia.ca, xiao.huang@concordia.ca}

% \title{An Integrated Matching Framework for Electric Vehicle  Ride-Hailing Services under Rider Demand and Charging Station Selection Uncertainty}

\title{Coordinating Guidance, Matching, and Charging Station Selection for Electric Vehicle Ride-Hailing Services through Data-Driven Stochastic Optimization}

% \affiliation[label1]{organization={},
% \affiliation[label2]{organization={},

\author[inst1]{Xiaoming Li}
\author[inst1]{Chun Wang}

\affiliation[inst1]{organization={Concordia Institute for Information Systems
Engineering (CIISE), Concordia University},%Department and Organization
            addressline={1455 Boulevard de Maisonneuve O}, 
            city={Montreal},
            postcode={H3G 1M8}, 
            state={Quebec},
            country={Canada}}

\author[inst2]{Xiao Huang}
\affiliation[inst2]{organization={John Molson School of Business (JMSB), Concordia University},%Department and Organization
            addressline={1455 Boulevard de Maisonneuve O}, 
            city={Montreal},
            postcode={H3G 1M8}, 
            state={Quebec},
            country={Canada}}

\begin{abstract}

Electric vehicles (EVs) play a pivotal role in sustainable ride-hailing services primarily due to their potential in reducing carbon emissions and enhancing environmental protection. Despite their significance, current research in the realm of EV batched matching frequently overlooks critical aspects such as rider demand uncertainty and charging station (CS) selection, leading to inefficiencies like decreased matching rates and prolonged waiting times for both riders and EV drivers. To fill the research gap, we propose a data-driven optimization framework that incorporates two inter-connected stochastic optimization models to address the challenges. The first model aims to relocate the idle EVs under satisfied conditions to the designated regions based on the probabilistic rider demand forecasting result before the real rider demand is revealed. Taking the solutions of the first model as the input, the second model optimizes the batched matching results by minimizing the rider's average waiting time and EV charging waiting time at CS. This integrated framework not only elevates the matching rate through the incorporation of rider demand uncertainties in the guidance module but also substantially curtails both rider and EV charging waiting times by synergizing guidance with CS selection choices. Empirical validation of our framework was conducted through an extensive case study in New York City, utilizing real-world data sets. The validation results demonstrate that the proposed data-driven optimization framework outperforms the benchmark models in terms of the proposed evaluation metrics. Most importantly, when deploying our framework, the charging waiting time of the EVs with low SOC can be reduced up to 73.6\% compared to the benchmark model without CS selection.

\end{abstract}

\begin{keyword}
EV ride-hailing \sep stochastic optimization \sep proactive guidance \sep batched matching \sep charging station selection

\end{keyword}

\end{frontmatter}


\section{Introduction}
With the rapid development of urbanization and population growth, ride-hailing services have undergone a significant transformation in shaping the urban
transportation paradigm. In the pursuit of sustainable mobility and the mitigation of environmental pollution, there has been a discernible shift toward the adoption of renewable energy vehicles, exemplified by the prevalence of electric vehicles (EVs) supplanting traditional Internal Combustion Engine Vehicles (ICEVs). Prominent Transportation Network Companies (TNCs), such as Lyft and Uber, have manifested a discernible trend in transitioning their ride-hailing services to EVs, as highlighted in the study by Jenn et al.~\cite{jenn2020emissions}.

EV plays a more pivotal role in shared mobility, especially within ride-hailing systems~\cite{pavlenko2019does}. The empirical findings reported in~\cite{jenn2020emissions} reveal that over the past decade, more than 50 million customers have availed themselves of ride-hailing services, encompassing a staggering 5.5 billion rides, thereby engendering a substantial volume of carbon emissions. The statistics from ~\cite{bauer2020demand, bauer2021leveraging} show that the growth of transportation emissions could easily offset the carbon reduction in other sections without rapid decarbonization. In this context, the adoption of EVs emerges as an efficacious strategy for ameliorating carbon emissions and mitigating air pollution.

Nevertheless, the transition from ICEVs to EVs is not devoid of challenges. Principally, the limited range of EVs will prevent the vehicle from being used in the same manner as ICEVs. Effective allocation of electric vehicles to riders necessitates judicious consideration of the state of charge (SoC) limitations. Furthermore, the proclivity of EVs to necessitate more frequent charging during operational hours introduces an exigency, as prolonged waiting times for charging can engender unavailability of EVs, thereby impinging upon the seamless provision of EV ride-hailing services. The intrinsic limitations of electric vehicles, characterized by reduced driving range, extended refueling durations compared to ICEVs, and constricted charging infrastructure, as elucidated in~\cite{pettit2019increasing}, pose formidable challenges to their widespread integration in the shared mobility landscape. Finally, the long waiting period for recharging will, in turn, impact the supply dynamics, necessitating the EV ride-hailing system to implement EV relocation operations.

There exist three inter-related key optimization challenges within the domain of EV ride-hailing services, namely, the \textit{EV relocation} (alternatively termed repositioning, rebalancing, and guidance in the literature) problem, the \textit{EV matching} (also known as order dispatching and vehicle assignment in the literature) problem, and \textit{EV charging station selection} (CS selection) problem. The relocation operation emerges as a pivotal solution for regulating the supply-demand equilibrium, exerting a considerable influence on the overall quality of ride-hailing services. Existing literature predominantly characterizes relocation as vehicle reactive repositioning and guidance as vehicle proactive repositioning, acknowledged as effective strategies for addressing imbalances in the supply-demand ratio~\cite{li2020data}. However, prevailing research predominantly operates within deterministic and static scenarios, often overlooking the inherent uncertainty associated with rider demand in the context of EV ride-hailing services. As pointed by~\cite{he2023data}, neglecting rider demand uncertainty detrimentally impacts the system-level performance of EV ride-hailing services. Additionally, extensive literature has been conducted to study the matching issue between EVs and riders~\cite{kullman2022dynamic, wang2023dynamic, chen2023electric, dong2022dynamic, gao2023stochastic}. Nonetheless, these studies often lack integration with EV relocation optimization, potentially resulting in low matching rates and prolonged rider waiting times. Lastly, while a substantial body of literature has explored the EV charging scheduling problem, the matter of CS selection within the context of EV ride-hailing services remains an area warranting more comprehensive investigation.

It is pertinent to highlight the intricate inter-dependence among the EV matching, EV relocation, and CS selection problems within the context of EV ride-hailing services. Addressing only one of these challenges is inherently incomplete, as the synergies between EV relocation, EV matching, and CS selection should be simultaneously considered. Notably, the outcomes of solving the charging problem, encompassing EV relocation and charging schedules, exert a consequential influence on the subsequent EV matching problem. Conversely, the results of the matching process, along with the updated vehicle routes derived from solving this problem, reciprocally impact the CS selection problem. Moreover, the surge in charging requests contributes to peak charging demand, thereby affecting EV availability and presenting challenges not only to EV relocation operations but also to the power grid. Given the high correlation among EV relocation, matching, and CS selection, a comprehensive approach is imperative to jointly address these three optimization issues in the context of EV ride-hailing services. While some works have endeavored to study the coordination of two or three issues in ride-hailing services~\cite{yang2023fleet, yi2021framework, zhan2022simulation, yu2023coordinating, dean2022synergies, guo2021robust}, these efforts often neglect considerations of uncertainty or the intricate interactions between decisions. This oversight renders such studies less practical for application in EV ride-hailing systems. Furthermore, uncertainty can be quantified or forecasted through statistical or machine learning approaches by leveraging historical data, which is frequently overlooked in existing literature.

The major contribution of this study is summarized as follows. We propose an integrated data-driven stochastic optimization framework that comprehensively addresses the challenges of idle EV proactive guidance under rider demand uncertainty and the EV matching problem under the uncertainty of CS selection. The framework comprises two distinct stochastic optimization models that collectively form an end-to-end solution to the identified problem. The first model, dedicated to proactive guidance strategy, systematically accounts for rider demand uncertainty, thereby generating a robust guidance solution. This solution effectively mitigates the supply-demand disparity within ride-hailing regions, ensures favorable matching rides with sufficient battery energy—and minimizes the idle travel distance of EVs. Subsequently, leveraging the optimal solution obtained from the first stochastic optimization model, the second model addresses batched matching solutions while accommodating CS selection uncertainty. This dual-model approach is adept at concurrently reducing rider waiting times and minimizing CS selection costs, which encompass EV travel distance and the anticipated waiting time of EVs at charging stations.

The rest of this paper is organized as follows: Section~\ref{section:literature-review} provides a comprehensive review of related work, Section~\ref{section:problem-statement} identifies the research gap and outlines the problem, Section~\ref{section:apporach} expounds upon the proposed data-driven stochastic optimization framework, Section~\ref{section:numerical-experiment} presents numerical validation, and finally, Section~\ref{section:conclusion} summarizes the contribution and discusses future research directions.

\section{Related Work}\label{section:literature-review}

In this section, we provide a comprehensive overview of the relocation, matching, and CS selection problems, along with their amalgamation in the context of ride-hailing services. The review initiates with an examination of pertinent studies involving ICEVs and subsequently transitions to investigations pertaining to EVs.

The relocation challenges associated with ICEVs have been extensively explored in the literature. Representative instances encompass relocation strategies such as Fluid model-based relocation~\cite{braverman2019empty}, queue theory-based relocation~\cite{zhang2016control}, Model Prediction Control (MPC)-based relocation~\cite{iglesias2018data}, and stochastic MPC-based relocation~\cite{tsao2018stochastic}. It is noteworthy that the majority of relocation strategies discussed in earlier works adhere to a reactive relocation approach, implying that the relocation processes are executed subsequent to the realization of rider demands in ride-hailing regions. Furthermore, ride-hailing matching strategies involving ICEVs, analyzed from an optimization perspective, encompass one-to-one matching~\cite{lowalekar2018online, zhang2017taxi} and one-to-many matching~\cite{li2021optimizing, meshkani2022generalized}. Subsequently, we delineate the decision-making problem studies focusing on EVs.

\subsection{Relocation in EV Ride-Hailing Services}
In contrast to relocation strategies employed for ICEVs, EVs necessitate careful consideration of battery levels during the relocation process. Zhao et al.\cite{zhao2018integrated} conducted a study on EV rebalancing and staff relocation within a one-way car-sharing system. The problem was formulated as a mixed-integer linear programming model aimed at minimizing EV rebalancing and staff relocation costs. The model underwent decomposition and resolution through the application of the Lagrangian relaxation technique. Gambella et al.\cite{gambella2018optimizing} addressed the EV relocation problem in a station-based, one-way trip car-sharing system. Formulating the reactive relocation problem as a mixed-integer programming model, the authors employed a rolling horizon approach for the solution. Qin et al.\cite{qin2022branch} delved into the EV relocation problem within one-way car-sharing systems, formulating it as a set-packing model and employing a customized branch-and-cut-and-price algorithm for resolution. Cai et al.\cite{cai2022hybrid} investigated the EV relocation problem with the objective of determining routing schedules for staff to rebalance the spatial-temporal distribution of EVs, optimizing for the highest total profit of the relocated requests. The proposed solution utilized a hybrid adaptive large neighborhood search and tabu search algorithm. Kim et al.\cite{kim2022idle} focused on shared autonomous EVs, employing deep learning models for predictive rider demand. The deep learning model's predictions were input for an EV relocation model. Huang et al.\cite{huang2020vehicle} undertook a comparative study between an operator-based model, where relocation is managed by the crew, and a user-based model, where price incentives are applied to balance supply and demand within a one-way station-based EV car-sharing system. Lin et al.~\cite{lin2021vehicle} investigated the EV deployment and relocation problem within a one-way station-based and free-floating car-sharing system under demand uncertainty, formulating the problem through a two-stage stochastic optimization model.

\subsection{Matching in EV Ride-Hailing Services}
Several studies have delved into the EV matching problem, exploring both one-to-one and one-to-many matching scenarios. These investigations encompass diverse methodologies, including deep reinforcement learning~\cite{kullman2022dynamic}, graph-based matching~\cite{wang2023dynamic}, a probabilistic matching approach~\cite{chen2023electric}, and models based on Markov decision processes~\cite{dong2022dynamic,gao2023stochastic}.

Kullman et al.\cite{kullman2022dynamic} addressed the EV ride-hailing service under central operator control, aiming to maximize platform profit through the judicious assignment of EVs to rider requests. The study employed a deep reinforcement learning algorithm utilizing the Q-value approximation technique for problem resolution. Wang et al.\cite{wang2023dynamic} proposed a two-stage graph-based many-to-one ride-matching algorithm tailored for shared autonomous EV systems. The first stage involved a bipartite graph-based one-to-one matching model with a pre-matching algorithm, while the second stage transformed the problem into a maximum weight matching scenario, establishing a vehicle-to-vehicle matching model. Chen et al.\cite{chen2023electric} introduced a probabilistic matching approach within EV ride-hailing services, aiming to maximize the platform's expected profit. The decision-making framework incorporated stochastic driver behaviors. Dong et al.\cite{dong2022dynamic} devised dynamic EV allocation policies for matching EVs to riders through a Markov decision process (MDP) model. The MDP parameters were derived from a semi-open queuing network with multiple synchronization stations, and the model sought to minimize rider waiting costs and lost demand, yielding a heuristic EV assignment policy. In a similar vein, Gao et al.~\cite{gao2023stochastic} explored a stochastic share-a-ride problem using EVs, considering a hybrid matching mode for customers and parcels. The problem was formulated as an MDP with the objective of maximizing overall expected rewards.

\subsection{Charging Station Selection in EV Ride-Hailing Services}
Several studies have also addressed the CS selection problem within the domain of EV ride-hailing services~\footnote{This section focuses on the impacts of charging management on ride-hailing services, excluding discussions related to EV charging scheduling or power grid considerations.}. Cao et al.\cite{cao2016ev} delved into the charging station selection problem, considering the trip duration and mobility uncertainties. The investigation incorporated considerations of EV charging reservations and parking durations at charging stations, with the overarching goal of minimizing drivers' trip durations through intermediate charging at selected stations. Yang et al.\cite{yang2021dynamic} examined the EV fast-charging station selection problem using a stochastic dynamic simulation modeling framework. The user behavior within the framework was formulated through a multi-nominal logit choice model, while the supply-demand dynamics were captured by a multi-server queuing model. Tian et al.~\cite{tian2016real} introduced a real-time Charging Station (CS) selection recommendation system tailored for EV ride-hailing services, leveraging large-scale GPS data mining. The real-time GPS trajectories facilitated CS selection based on recharging start times. It is noteworthy that the chosen CSs significantly impact subsequent relocation within the operational window. Hence, it is rational to concurrently consider the coordination between these two decision-making problems.

\subsection{Integration of Multiple Decision-Making Problems}
Several studies have focused on the integration of multiple decision-making problems within EV ride-hailing services. Yang et al.\cite{yang2023fleet} investigated fleet sizing and charging infrastructure design, employing three distinct mathematical models to address vehicle matching, relocation, and charging assignment. A Bayesian optimization framework was applied to optimize fleet size and charging station configuration jointly. Yi et al.\cite{yi2021framework} examined EV dispatching and its interaction with charging facilities in ride-hailing services, proposing optimization models and heuristic rules to simulate the EV ride-hailing framework. Zhan et al.\cite{zhan2022simulation} introduced a simulation-optimization framework for dynamic EV ride-hailing services, combining the EV ride-hailing matching problem and EV charging problem through two optimization models—one formulated as a mixed-integer linear programming model and the other as a mixed-integer nonlinear programming model. Yu et al.\cite{yu2023coordinating} explored the coordination of matching, rebalancing, and charging optimization problems for EV ride-hailing services under hybrid rider request modes. Dean et al.\cite{dean2022synergies} investigated synergies between EV repositioning and charging strategies, proposing an optimization framework to minimize waiting times and the fleet's percentage of empty travel. Guo et al.~\cite{guo2021robust} explored ride-hailing matching and vehicle rebalancing problems using a robust optimization model that considers uncertainty in rider demand, aiming to minimize overall costs in terms of total vehicle miles traveled and the number of unsatisfied requests under worst-case demand scenarios.

Existing works have often considered only a subset of decision-making problems or overlooked uncertainty. As discussed, these three decision-making problems exhibit high correlation, and the integration of operations can yield significant benefits for EV ride-hailing services. Moreover, the stochastic and dynamic nature of EV ride-hailing services suggests that incorporating uncertainty into the platform operations can be advantageous. The proposed data-driven optimization framework is elaborated upon in the subsequent section.

\begin{table*}[htbp]
\caption{Summary of existing work that integrates multiple decision-making problems}
\renewcommand{\arraystretch}{1.3}  %
\centering
\small % This sets the font size to small,  \tiny, \scriptsize, \footnotesize,\small, \normalsize, \large, \Large (the 'L' is in uppercase), \LARGE (all uppercase), \huge, \Huge (the 'H' is in uppercase)
\setlength{\tabcolsep}{2.1mm}{
% \scriptsize     % font size in table
\begin{tabular}{cccccc}
\toprule[1pt]
Ref. & Vehicle Relocation & Trip Matching & CS Selection & Uncertainty Involvement & Learning-Based\\ 
\hline

{\cite{yang2023fleet}} & {Reactive} & {$\checkmark$} & {$\checkmark$}  & {$\usym{2717}$} & {$\usym{2717}$} \\
{\cite{yi2021framework}} & {Reactive} & {$\checkmark$} & {$\checkmark$}  & {$\usym{2717}$} & {$\usym{2717}$} \\ 
{\cite{zhan2022simulation}} & {$\usym{2717}$} & {$\checkmark$} & {$\checkmark$}  & {$\usym{2717}$} & {$\usym{2717}$} \\ 
{\cite{yu2023coordinating}} & {Reactive} & {$\checkmark$} & {$\usym{2717}$}  & {$\usym{2717}$} & {$\usym{2717}$} \\
{\cite{dean2022synergies}} & {Reactive} & {$\usym{2717}$} & {$\checkmark$}  & {$\usym{2717}$} & {$\usym{2717}$} \\
{\cite{guo2021robust}} & {Reactive} & {$\checkmark$} & {$\usym{2717}$}  & {Rider Demand} & {$\usym{2717}$} \\
\hline
{This work} & {Proactive} & {$\checkmark$} & {$\checkmark$} & {Rider Demand, CS Selection} & {Probabilistic Forecasting}\\ 
\bottomrule[1pt]
\end{tabular}}
\label{table:literature}
\end{table*}

\section{The EV Ride-Hailing Problem Statement}\label{section:problem-statement}

In this section, a comprehensive exploration of the studied EV ride-hailing platform is presented. Initially, we provide a concise delineation of the pivotal components constituting the EV ride-hailing platform. Subsequently, we explain how the studied EV ride-hailing platform operates systematically.

\subsection{System Overview}
We consider a fairly general EV ride-hailing service platform which provides ride-hailing services in a set of regions. The platform involves three key entities: a fleet of EVs with associated drivers~\footnote{For simplicity, the terms EV and driver will be used interchangeably henceforth.}, a set of riders, and a central operator responsible for platform operations. Equipped with sensors and wireless connectivity, the central operator can monitor the real-time position and status of EVs, including vacancy status and State of Charge (SoC). Simultaneously, riders can submit service requests via mobile devices or the Internet, with the central platform managing and aggregating data for subsequent decision-making. Furthermore, the central platform maintains information regarding CSs in proximity to the ride-hailing regions. The systematic integration of these data sources forms the basis for addressing decision-making problems.

\subsection{Problem Description}

\begin{figure*}
    \centering
    \includegraphics[scale=0.7]{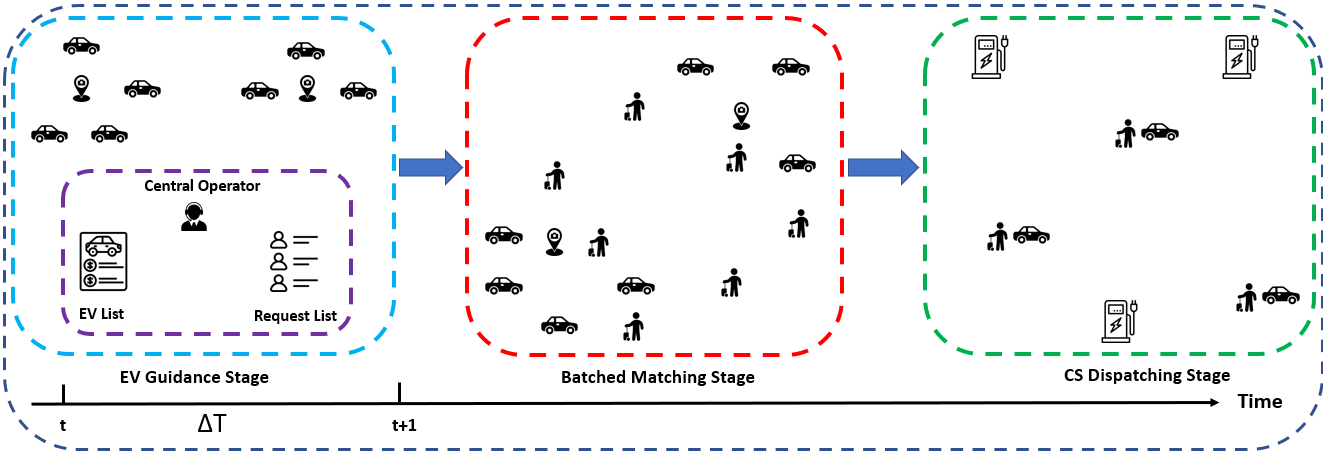}
    \caption{The diagram of the proposed approach for EV ride-hailing system}
    \label{fig:system-overview}
\end{figure*}

The operational framework follows a periodic rhythm defined by a time interval termed a \textit{batching window}. The daily timeline is discretized into batching windows of fixed duration, denoted as $\Delta T$ (e.g., 10 minutes). As depicted in Fig.\ref{fig:system-overview}, at the commencement of each batching window (time $t$), the platform orchestrates the movement of vacant EVs, possessing ample SoC, to predetermined ride-hailing regions with specified points of interest (POIs, e.g., landmarks) to meet rider demand. Concurrently, the platform aggregates ride requests from riders. It is noteworthy that idle EVs are positioned in regions ahead of the actual rider demand realization. Ideally, the guided EVs (supply) would precisely match the rider requests (demand). However, inherent biases in supply and demand necessitate reliance on a probabilistic forecasting approach for demand. Several works, such as\cite{chen2023probabilistic}, have delved into probabilistic demand forecasting within the transportation domain. 

After the EV guidance operation and rider request collection, the central platform embarks on optimizing the one-to-one matching problem at time $t+1$. Unlike the conventional ride-hailing matching using ICEVs, where the primary goal is to minimize rider waiting time, our work incorporates the decision-making process for CS selection during matching. Following completion of the trip (rider drop-off), EVs then traverse to the assigned CS. The entire EV ride-hailing platform operates iteratively in a rolling horizon manner.

\section{The Integrated Solution Framework}\label{section:apporach}

In this section, we elaborate on the proposed data-driven optimization framework, which encompasses duo stochastic optimization models. The first model is to optimize the idle EV proactive guidance under the rider demand uncertainty; the output of the stage is the optimal guidance solution. Taking the solutions as the input, the second model is to match the drivers and riders while considering the CS selection uncertainty. The notations are shown in Table~\ref{table:notation}.

% \begin{comment}
\begin{table}[htbp]
\caption{Notations in Mathematical Optimization Models}
\renewcommand{\arraystretch}{1.3}  %
\small % This sets the font size to small,  \tiny, \scriptsize, \footnotesize,\small, \normalsize, \large, \Large (the 'L' is in uppercase), \LARGE (all uppercase), \huge, \Huge (the 'H' is in uppercase)
\centering
\setlength{\tabcolsep}{.9mm}{
% \scriptsize     % font size in table
\begin{tabular}{ll}
\toprule[1pt]
\textbf{Indices} & \textbf{Description} \\ 
\hline
% {$i$} & {The index of taxi} \\
{$i$} & {The index of ride-hailing regions} \\
{$j$} & {The index of EVs}  \\
{$k$} & {The index of riders} \\ 
{$u$} & {The index of charging stations} \\ 
% {$s$} & {The index of demand scenarios} \\ 
{$t$} & {The index of batching windows} \\ 
\hline
{\textbf{Sets}} & {\textbf{Description}} \\ \hline
% {$\mathscr{M}$} &{A set of vacant taxis indexed by $i$} \\
{$\mathcal{A}$} & {A set of ride-hailing regions index by $i$}  \\
{$\mathcal{D}$} & {A set of idle EV index by $j$} \\
{$\mathcal{R}$} & {A set of riders indexed by $k$} \\
{$\mathcal{U}$} & {A set of charging stations indexed by $u$} \\
% {$\mathcal{S}$} & {A set of demand scenarios indexed by $s$} \\
{$\mathcal{T}$} & {A set of batching windows indexed by $t$} \\\hline 
{\textbf{Params.}} & {\textbf{Description}} \\ \hline
{$\hat{d}_{i}^{t}$} & {The sample (outcome) of the forecasting rider demand in the region $i$ at time $t$}  \\
% {} & {at time slot $t$ under scenario $s$}\\
{$\widehat{D}_{i}^{t}$} & {The random variable of the forecasting rider demand in region $i$ at time $t$}  \\
{$f_{\widehat{D}_{i}^{t}}(\cdot)$} & {The probability density function (PDF) of the random variable above}  \\
{$\mu(\cdot)$} & {A probability measure of the forecasting demand} \\
{$g_{i, j}$} & {The guidance distance between region $i$ and EV $j$} \\
{$\mathcal{C}^{t}$  } & {The number of idle EVs at time $t$} \\
{$\alpha$} & {The coefficient of EV's idle distance cost } \\
{$\beta_{1}, \beta_{2}$} & {The over-supply and under-supply cost coefficients, respectively} \\
{$\gamma$} & {The EVs travel speed} \\
{$\lambda$} & {The minimum SoC coefficient} \\
{$\omega_{j}$} & {Battery consumption rate of EV $j$} \\
{$\Delta T$} & {The fixed duration of batching window} \\
{$ST$} & {The starting time of batching window} \\
{$H$} & {A large positive number} \\
{$\theta_{1}$} & {The coefficient of CS selection cost} \\
{$\theta_{2}$} & {The coefficient of rider's waiting time} \\
{$\pi_{1}$} & {The coefficient of travel distance to CS} \\
{$\pi_{2}$} & {The coefficient of waiting time at CS} \\
{$o(k), w(k)$} & {The origin and destination of rider $k$} \\
{$c(\cdot)$} & {The coordinate of POI or CS} \\
{$W_{j, k}^{1}$} & {The charging cost of EV $j$ severs rider $k$} \\
{$W_{j, k}^{2}$} & {The waiting time of rider $k$ pick up by EV $j$} \\
{$req_{k}$} & {The ride request time of rider $k$} \\
{$trip_{i}^{avg}$} & {The average trip distance at region $i$} \\
{$trip_{k}$} & {The trip distance of rider $i$} \\
{$dist_{\cdot(\cdot), \cdot(\cdot)}$} & {The distance between two coordinates} \\
{$ldt_{k}$} & {The latest departure time of rider $k$} \\
{$WT_{u}^{t}$} & {The random variable of waiting time in charging station $u$ at time $t$} \\
{$SoC_{j}$} & {The State of Charge of EV $j$ before guidance} \\
{$SoC_{j}^{\prime}$} & {The State of Charge of EV $j$ after guidance} \\
\hline
{\textbf{Vars.}} &{\textbf{Description}}\\\hline
{$x_{i, j}^{t}$} & {1 if EV $j$ is guided to region $i$ at time $t$, 0 otherwise}\\
% {} & {  0 otherwise} \\ 
{$y_{j, k}$} & {1 if EV $j$ is matched to rider $k$, 0 otherwise }\\
\bottomrule[1pt]
\end{tabular}}
\label{table:notation}
\end{table}
% \end{comment}

\subsection{The Idle EV Guidance Model}
Let $\mathcal{A}$, $\mathcal{D}$, and $\mathcal{R}$ denote a set of EV ride-hailing regions in a given area, a set of available EV drivers, and a set of riders, respectively~\footnote{In fact, $\mathcal{D}^{t}$ and $\mathcal{R}^{t}$ denote the EV drivers and riders at the current time $t$. Since the superscript $t$ does not change the problem characteristics, we omit $t$ for simplicity.}. During each batching window, a number of idle EVs are guided to the different regions to satisfy future rider requests. Since the proactive guidance requires the idle EVs to move to the designed regions before the real demands are realized, it is imperative to accurately forecast the rider demand for the current batching window. Notice that the rider demand in the ride-hailing system follows a time-series pattern. Therefore, this can be archived by time-series forecasting approaches. Two types of forecasting approaches can be applied to rider demand prediction, namely, point prediction and probabilistic prediction. The former returns a scalar as the outcome, while the latter returns a distribution as the result. Due to the stochastic optimization model characteristics, the rider demand distribution is required in our problem setting. Notice that the rider demand distribution can be predicted by the probabilistic forecasting approach, we adopt our previous work as the learning component~\cite{li2023xrmdn}.

Meanwhile, rider $k \in \mathcal{R}$ sends a request to the central platform. The request format is characterized as a 5-tuple (Request-ID, o(k), w(k), $req_{k}$, $ldt_{k}$), where Request-ID is a unique string, o(k) and w(k) denote the origin and destination of rider $k$, respectively, $req_{k}$ represents the time-stamp of the request submission by rider $k$, and $ldt_{k}$ means the latest departure time of rider $k$. In addition, we define a binary decision variable $x_{i, j}^{t} = 1$, if driver $j$ is guided to region $i$, and $x_{i, j}^{t} = 0$ otherwise.

The first objective of the proactive guidance model is to minimize the overall idle moving costs of the guided EV drivers. It is reasonable to assume the cost is proportional to the EV idle moving distance. Let $g_{i,j}$ be the idle driving distance from the driver $j$'s current location to the POI of region $i$, the total costs that are associated with the EV idle guidance distances can be denoted by:
\begin{equation}\label{SO-1-guidance}
\begin{aligned}
 \sum_{i \in \mathcal{A}} \sum_{j \in \mathcal{D}} \alpha_{j} \cdot g_{i, j} \cdot x_{i, j}^{t}
\end{aligned} 
\end{equation}

where $\alpha_{j}$ represents the idle moving cost per mile of EV driver $j$. In addition,  the costs incurred by the imbalance of supply (i.e., the number of guided EVs) and demand (i.e., the number of predicted rider requests) of the ride-hailing regions, including over-supply costs and under-supply cost, can be represented as:

\begin{equation}\label{SO-1-over-supply}
\begin{aligned}
\beta_{1} \sum_{i \in \mathcal{A}} \int_{\Omega} \max \left\{0, \sum_{j \in \mathcal{D}} x_{i, j}^{t} - \hat{d}_{i}^{t} \right\} f_{\widehat{D}_{i}^{t}} (\hat{d}_{i}^{t}) \mu(\hat{d}_{i}^{t}) 
\end{aligned} 
\end{equation}

\begin{equation}\label{SO-1-under-supply}
\begin{aligned}
\beta_{2} \sum_{i \in \mathcal{A}} \int_{\Omega} \max \left\{0, \hat{d}_{i}^{t} - \sum_{j \in \mathcal{D}} x_{i, j}^{t} \right\} f_{\widehat{D}_{i}^{t}} (\hat{d}_{i}^{t}) \mu(\hat{d}_{i}^{t})
\end{aligned} 
\end{equation}

where $\hat{d}_{i}^{t}$ and $f_{\widehat{D}_{i}^{t}}$ denote the predicted rider demand scenario (sampled outcome) in region $i$ at time $t$ and the probability density function (PDF) of the random variable of the predicted rider demand $\widehat{D}_{i}^{t}$ in region $i$ at time $t$, respectively. Eventually, combining the two terms above, the objective function of the idle EV guidance model can be denoted as:

\begin{equation}\label{SO-1-obj}
\begin{aligned}
\sum_{i \in \mathcal{A}} \sum_{j \in \mathcal{D}} \alpha_{j} \cdot g_{i, j} \cdot x_{i, j}^{t} + \beta_{1} \sum_{i \in \mathcal{A}} \int_{\Omega} \max \left\{0, \sum_{j \in \mathcal{D}} x_{i, j}^{t} - \hat{d}_{i}^{t} \right\} f_{\widehat{D}_{i}^{t}} (\hat{d}_{i}^{t}) \mu(\hat{d}_{i}^{t}) + \beta_{2} \sum_{i \in \mathcal{A}} \int_{\Omega} \max \left\{0, \hat{d}_{i}^{t} - \sum_{j \in \mathcal{D}} x_{i, j}^{t} \right\} f_{\widehat{D}_{i}^{t}} (\hat{d}_{i}^{t}) \mu(\hat{d}_{i}^{t})
\end{aligned} 
\end{equation}

where $\beta_{1}$ and $\beta_{2}$ are positive numbers to control the weight of over-supply and under-supply imbalance costs, respectively. Notice that the integral terms in the objective function can be linearized via Monte Carlo family algorithms, such as sample average approximation~\cite{birge2011introduction}. Thus, the reformulated model can then be solved by off-the-shelf mathematical solvers. After discussing the objective function, a few constraints must be satisfied. During the current batching window, each EV can be guided to one ride-haling region at most, which can be captured by:

\begin{equation}\label{constraint:c3}
\begin{aligned}
\sum_{i \in \mathcal{A}} x_{i, j}^{t} \leq 1, \quad \forall j \in \mathcal{D}, \quad \forall t \in \mathcal{T}
\end{aligned}
\end{equation}

As discussed in Section~\ref{section:problem-statement}, the central operator can monitor the status of EVs. The value can be estimated by the average of EVs from historical data. The fleet size is determined by the following constraint:

\begin{equation}\label{constraint:c4}
\begin{aligned}
\sum_{i \in \mathcal{A}} \sum_{j \in \mathcal{D}} x_{i, j}^{t} \leq \mathcal{C}^{t}, \quad \forall t \in \mathcal{T}
\end{aligned}
\end{equation}

In addition, EV $j$, if it can be guided to region $i$, must be able to reach the corresponding POI of the region within the batching window (namely, $\Delta T$). This can be denoted by the following constraint: 

\begin{equation}\label{constraint:c5}
\begin{aligned}
g_{i, j} / \gamma \leq \Delta T + H (1 - x_{i, j}^{t}),  \quad \forall i \in \mathcal{A}, \quad \forall j \in \mathcal{D},  \quad \forall t \in \mathcal{T}
\end{aligned}
\end{equation}

where $\gamma$ is the EV travel speed assumed to be a constant value, and $H$ is a large positive number to represent the logic "if"~\cite{rardin1998optimization}. Further, if EV $j$ can be guided to region $i$, it must guarantee that the SoC is sufficient to travel to the region $i$ to serve the rider trip. The constraints can be represented by:

\begin{equation}\label{constraint:c6}
\begin{aligned}
\omega_{j} \left(g_{i, j} + trip_{i}^{avg}\right) + \lambda \cdot SoC_{j} \leq SoC_{j} + H \left(1 - x_{i, j}^{t} \right), \quad \forall j \in \mathcal{D}, \quad \forall t \in \mathcal{T}
\end{aligned}
\end{equation}

where $trip_{i}^{avg}$ denotes the rider's average trip distance that starts from region $i$. The parameters $\omega_{j}$ and $\lambda$ denote the battery consumption per kilometer and minimum proportion of SoC (reserved battery for the future CS selection), respectively. Notice that the parameter $trip_{i}^{avg}$ is unknown during the guidance operation, however, the value of this parameter can be learnt or estimated from historical data. In this sense, statistical learning or machine learning approaches would be promising methods to achieve that goal.

Given the decision variables, parameters, objective function, and constraints discussed above, the EV proactive guidance problem is formulated as a stochastic integer linear optimization model, shown as follows.

\begin{equation}\label{model-1}
\begin{aligned}
& \operatorname{minimize} \quad  (\ref{SO-1-obj})\\
& \text {s.t. (\ref{constraint:c3}), (\ref{constraint:c4}), (\ref{constraint:c5}), (\ref{constraint:c6}) } \\
& x_{i, j}^{t} \in \{0, 1\}, \quad \forall i \in \mathcal{A}, \quad \forall j \in \mathcal{D}, \quad \forall t \in \mathcal{T}
\end{aligned}
\end{equation}

\subsection{The EV Batched Matching Model with Charging Station Selection}
After the batching window ends, the batching matching component starts to operate. Given the output solutions from the EV idle guidance model as the input of the matching model, the EVs and riders are matched accordingly. One major objective in the batched matching model is to minimize the CS selection, which involves the EV travel distance to the CS and the expected waiting time at the CS. If rider $k$ is matched with EV driver $j$, let $\mathcal{U}_{k}$ be the set of CSs that surrounds the destination of rider $k$ (e.g., CSs scattered one mile around the rider's destination), the EV $j$'s travel distance between the rider $k$'s destination and the CS $u \in \mathcal{U}_{k}$ is denoted by:

\begin{equation}
\begin{aligned}
\pi_{j}^{1} \cdot dist_{w(k), c(u)} 
\end{aligned}
\end{equation}

where $\pi_{j}^{1}$ is the EV's travel cost per kilometer. Another critical purpose of the CS selection is that we expect the EVs with low batteries to be sent to the CS, where the waiting time is short. This assumption makes sense that the EVs with low batteries can be recharged in a short time~\cite{jamshidi2021dynamic}. Therefore, it can be used for the EV ride-hailing service. This can be illustrated by Fig.~\ref{fig:CS-selection}. In this case, compared to the EV with high SoC (80\%), the EV with lower SoC (30\%) has the priority to be recharged in a shorter time, therefore, it is assigned to the rider whose destination is in the green rounded rectangle. This is because the minimum expected waiting time of CS (3 mins) in the green rounded rectangle is shorter than the one (8 mins) in the red rounded rectangle. Eventually, the costs of waiting time at the CS lead to the following term:

\begin{figure}
    \centering
    \includegraphics[scale=.7]{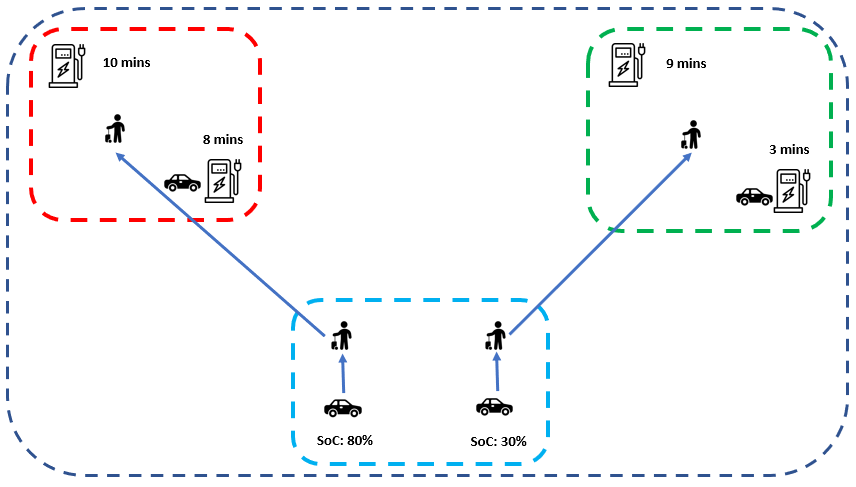}
    \caption{An example of CS selection during the batched matching optimization, EV with low SoC has the priority to match the rider whose destination deploys CS with short expected waiting time}
    \label{fig:CS-selection}
\end{figure}

\begin{equation}\label{SO-matching-obj:term-2}
\begin{aligned}
\pi_{j}^{2} \cdot \frac{\EX [WT_{u}^{t}]}{SoC_{j}^{\prime} - \omega_{j} \cdot trip_{k}} 
\end{aligned}
\end{equation}

where $\pi_{j}^{2}$ and $\omega$ denote the coefficient of waiting time cost at CS and the $SoC$ consumption per kilometer, respectively. Combing the two weighted terms together, the minimum costs of driver $j$ choosing the optimal CS selection that is located around the destination of rider $k$ can be defined as follows:

\begin{equation}\label{constraint:c10}
\begin{aligned}
W_{j, k}^{1} := \min_{u \in \mathcal{U}_{k}} \left\{\pi_{j}^{1} \cdot dist_{w(k), c(u)} + \pi_{j}^{2} \cdot \frac{\EX [WT_{u}^{t}]}{SoC_{j}^{\prime} - \omega_{j} \cdot trip_{k}} \right\}
\end{aligned}
\end{equation}

In addition, the rider's waiting time is a common criterion to measure the service quality of a ride-hailing system. In this work, rider $k$'s waiting time is defined as the difference between EV driver $j$'s arrival time (the time that EV $j$ pick up rider $k$) and rider $k$'s requested time $req_{k}$. Therefore, if driver $j$ is matched to rider $k$, the rider's waiting time can be represented as follows:

\begin{equation}\label{constraint:c11}
\begin{aligned}
W_{j, k}^{2} :=  ST + \Delta T - req_{k} + \frac{dist_{c(j), o(k)}}{\gamma}
\end{aligned}
\end{equation}

Unlike in the existing literature, where the matching objective considers rider waiting time only. Charging waiting time is also incorporated in this study, which means both riders and EV drivers can benefit from our batching matching model. In fact, this optimization model is from the perspective of a two-sided market. Essentially, it provides a trade-off solution for the EV ride-hailing systems. The equilibrium solutions enable the platform to be more practical. Let $y_{j,k}$ be the binary decision variable, $y_{j, k} = 1$, if EV driver $j$ is matched with rider $k$, and $y_{j, k} = 0$ otherwise. Eventually, the objective function~(\ref{SO-2-obj}) can be formulated as follows: 

\begin{equation}\label{SO-2-obj}
\begin{aligned}
\sum_{j \in \mathcal{D}} \sum_{k \in \mathcal{R}} (\theta_{1} W_{j, k}^{1} + \theta_{2} W_{j, k}^{2}) \cdot y_{j, k} + H \left(|\mathcal{R}| - \sum_{j \in \mathcal{D}} \sum_{k \in \mathcal{R}} y_{j, k} \right)
\end{aligned}
\end{equation}

where the second term in the objective function

\begin{equation}
\begin{aligned}
H \left(|\mathcal{R}| - \sum_{j \in \mathcal{D}} \sum_{k \in \mathcal{R}} y_{j, k} \right)
\end{aligned}
\end{equation}

denotes the penalty incurred by the unmatched requests, which is controlled by a big positive integer $H$. Besides the objective function, a few constraints must be satisfied in the batched matching model. Since we consider the one-to-one matching problem in this study, namely, each driver can serve at most one rider, meanwhile each rider can be assigned to one driver at most. This leads to the following pair of constraints: 

\begin{equation}\label{constraint:c12}
\begin{aligned}
\sum_{k \in \mathcal{R}} y_{j, k} \leq 1, \quad \forall j \in \mathcal{D}
\end{aligned}
\end{equation}

\begin{equation}\label{constraint:c13}
\begin{aligned}
\sum_{j \in \mathcal{D}} y_{j, k} \leq 1, \quad \forall k \in \mathcal{R}
\end{aligned}
\end{equation}

Notice that the constraints above imply that not all the drivers and riders are mandatory to be matched. In addition, if EV $j$ is matched with rider $k$, EV $j$ must be able to complete the rider $k$'s trip and travel to the furthest CS which is located around the rider $k$'s destination to recharge after dropping the rider off. This results in the following group of constraints.  

\begin{equation}\label{constraint:c14}
\begin{aligned}
\omega_{j} \cdot \left(trip_{k} + \max_{u \in \mathcal{U}_{k}} dist_{w(k), c(u)}\right) \leq SoC_{j}^{\prime} + H (1 - y_{j, k}),  \quad \forall j \in \mathcal{D}, \quad \forall k \in \mathcal{R}
\end{aligned}
\end{equation}

where $dist_{w(k), c(u)}$ denote the distance between rider $k$'s destination and location of CS $u$. Also, as illustrated in Fig.~\ref{fig:RWT}, if EV driver $j$ can be matched with rider $k$, EV $j$'s arrival time (the time driver $j$ picks up rider $k$) must be earlier than rider $k$'s latest departure time $ldt_{k}$. This condition can be ensured by the following constraints:

\begin{equation}\label{constraint:c15}
\begin{aligned}
ST + \Delta T + \frac{dist_{c(j), o(k)}}{\gamma} \leq ldt_{k} + H (1 - y_{j, k}),  \quad \forall j \in \mathcal{D}, \quad \forall k \in \mathcal{R}
\end{aligned}
\end{equation}

where $dist_{c(j), o(k)}$ denotes the distance between EV $j$'s starting location and rider $k$'s destination. Given the decision variables, parameters, objective function, and constraints discussed above, the EV batch matching problem is formulated as a stochastic integer linear optimization model, shown as follows.

\begin{equation}\label{model-2}
\begin{aligned}
& \operatorname{minimize} \quad  (\ref{SO-2-obj})\\
& \text {s.t. (\ref{constraint:c10}), (\ref{constraint:c11}), (\ref{constraint:c12}), (\ref{constraint:c13}), (\ref{constraint:c14}), (\ref{constraint:c15}) } \\
& y_{j, k} \in \{0, 1\}, \quad \forall j \in \mathcal{D}, \quad \forall k \in \mathcal{R}
\end{aligned}
\end{equation}

Since the EV ride-hailing platform operates in a rolling horizon manner, if there exist a few riders that are not matched in the current stage, in this study, we assume that they will wait for the next batching matching for simplicity. Although the matching model is formulated as an integer programming optimization problem, it can be effectively solved by off-the-shelf mathematical solvers or user-defined algorithms, such as the Hungarian algorithm~\cite{kuhn1955hungarian}.

\section{Numerical Experiment}\label{section:numerical-experiment}
In this section, we carry out a comprehensive study to validate our proposed framework. We begin by specifying the experiment setup, including model parameters and experiment environment. We then proceed to explain the structure of the dataset and the data processing that was carried out. Following this, we elaborate on the evaluation metrics and benchmark models adopted in the study. Finally, we compare the performance of our approach to the benchmark models in terms of the specified performance criteria.

\subsection{Experiment Setup}
We assume there are three types of EVs with different battery consumption rates, which is the same configuration in~\cite{cao2016ev}. In addition, the EVs are assumed to travel at a constant speed of 30 kpm. According to the report by the American Automobile Association in 2022~\cite{AAA_2022}, the costs of hybrid and electric vehicles are 64.61\textcent/mile and 60.32\textcent/mile, respectively. Based on that, the EV guidance cost $\alpha_{j}$ in the following experiment is assumed to follow a uniform distribution in the range of \$0.8 to \$1.1 per kilometer. Moreover, we expect more riders can be served by the EV ride-haling services, hence, the over-supply cost $\beta_1$ is less than the under-supply cost $\beta_2$, which are set to 5 and 10, respectively. Notice that we believe the charging waiting time is much more critical than the EV's travel distance to the designated CS in the CS selection costs. Therefore, the two weights $\pi_2$ and $\pi_1$ are set to 10 and 1, respectively. Eventually, the parameters in the optimization models are summarized in Table~\ref{table:model-parameter}.

\begin{table}[htbp]
\caption{Parameter settings in the optimization models}
\renewcommand{\arraystretch}{1.3}  %
\centering
\small % This sets the font size to small,  \tiny, \scriptsize, \footnotesize,\small, \normalsize, \large, \Large (the 'L' is in uppercase), \LARGE (all uppercase), \huge, \Huge (the 'H' is in uppercase)
\setlength{\tabcolsep}{5.5mm}{
% \scriptsize     % font size in table
\begin{tabular}{lll}
\toprule[1pt]

{\textbf{Parameters}} & {\textbf{Value}} & {\textbf{Description}} \\ \hline
{$\alpha_{j}$} & {\$0.8/km - \$1.1/km} & {Guidance cost coefficient} \\
{$\beta_{1}$} & {5} &  {The over-supply cost coefficient} \\
{$\beta_{2}$} & {10} &  {The under-supply cost coefficient} \\
{$\gamma$} & {30 kph} &  {EVs travel speed} \\
{$\lambda$} & {10\%} &  {Minimum SoC of EV} \\
{$\omega_{j}$} & {\{ 0.1171 kWh/km, 0.1751 kWh/km, 0.1863 kWh/km \}} &  {Battery consumption per kilometer} \\
{$\Delta T$} & {10 mins} &  {Duration of time window} \\
{$\theta_{1}$}  & 1 & {The weight of CS selection cost} \\
{$\theta_{2}$}  & 10 & {The weight of rider waiting time} \\
{$\pi_{1}$} & 1 & {The weight of distance cost} \\
{$\pi_{2}$} & 10 & {The weight of charging waiting time} \\
{$SoC_{j}$} & {$\mathcal{U}$ (0.2, 0.8)} & {The initial EV's State of Charge} \\
% {} & {} & {parameter setting, see~\cite{yang2021dynamic} } \\
\bottomrule[1pt]
\end{tabular}}
\label{table:model-parameter}
\end{table}

The mathematical models are coded in Python 3.11 and solved by Gurobi 11.0~\footnote{https://www.gurobi.com/downloads/gurobi-software/}. The experiments are run on a PC with Intel Core i7, 32 GB RAM, and Windows 11. The average computational time of the guidance model and matching model are around 17.6 s and 29.8 s, respectively. The result suggests that our framework can be solved in a reasonable time and can be deployed on real EV ride-hailing systems.

\subsection{Data Description and Data Process}

\begin{figure}[!htb]
\minipage{0.32\textwidth}
  \includegraphics[width=\linewidth]{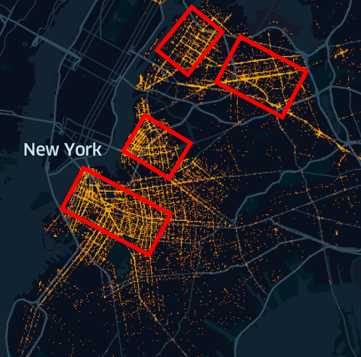}
  \caption{Heatmap of Rider Pick-Up Density}\label{fig:pick-up}
\endminipage\hfill
\minipage{0.32\textwidth}
  \includegraphics[width=\linewidth]{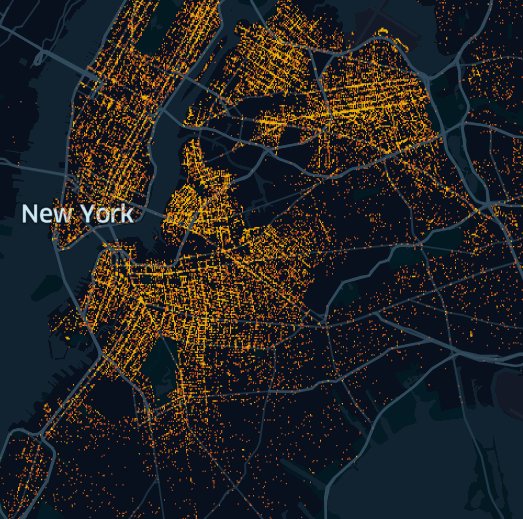}
  \caption{Heatmap of Rider Drop-Off Density}\label{fig:drop-off}
\endminipage\hfill
\minipage{0.32\textwidth}%
  \includegraphics[width=\linewidth]{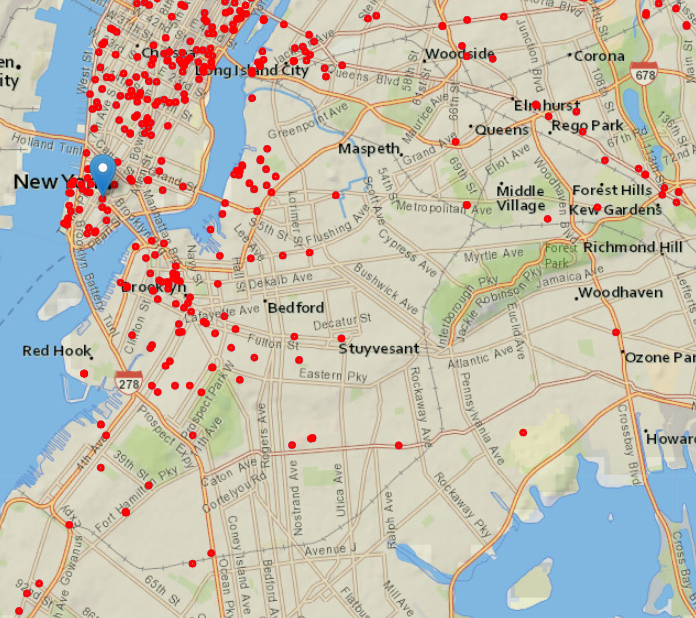}
  \caption{Charging Station Distribution}\label{fig:charging-station}
\endminipage
\end{figure}

Two real data sets - New York green taxi trip records~\footnote{The up-to-date data sets in parquet format (https://www.nyc.gov/site/tlc/about/tlc-trip-record-data.page) do not contain the coordinate information, which cannot be used in our problem setting since the coordination information is required. Instead, we leverage the archived data sets with the rider's pick-up and drop-off coordinates, which are available upon requests}, and charging stations in New York~\footnote{https://data.ny.gov/Energy-Environment/Electric-Vehicle-Charging-Stations-in-New-York/7rrd-248n/about\_data} are leveraged for the experiment validation. The taxi trip record data set in January 2016 contains 454,273 trip records and 20 attributes. We filtered a subset of attributes, including rider pickup \& drop-off datetime, pickup \& drop-off latitude \& longitude and applied them in the following experiment. In addition, a rider's trip distance is calculated by the Euclidean distance between the pickup and drop-off coordinates. Meanwhile, we randomly draw 5,000 samples of rider pickup and drop-off coordinates from the trip records and plot the heatmaps in Fig.~\ref{fig:pick-up} and Fig.~\ref{fig:drop-off}, respectively. Based on the visualization result, four regions (labeled by red rectangles) with high density are selected as the EV ride-hailing service regions. The fixed length of the batching window is set to 10 minutes, therefore, there are 144 batching windows in one day, and the trip records are aggregated every 10 minutes as the rider demand of the current batching window.

As a data-driven optimization framework, the guidance model in the framework requires predicted rider demand as the inputs, we apply AIRMA as the rider demand ($\hat{d}_{i}^{t}$) forecasting approach for the deterministic guidance model, and XRMDN~\cite{li2023xrmdn} as the probabilistic rider demand ($\widehat{D}_{i}^{t}$) forecasting approach for the stochastic guidance model. We leverage the trip records in the last three months of 2015 and derive the rider demands by aggregating the trip records every 10 minutes, which is introduced in the previous paragraph. The historical rider demand information is then used as the training sets for ARIMA and XRMDN.

In addition, there are 3,652 charging stations in the raw data set, we select 388 charging stations in the selected area as shown in Fig.~\ref{fig:charging-station}. The only features required in our experiment are the charging station coordinates. As the case study in the greater Seoul metropolitan area~\cite{ko2017determining} indicates, 4 km service distance can achieve over 90\% charging service coverage rate. Considering the trip data in the New York metropolitan area, we assume that when the EV completes the trip (drops off the matching rider), the driver will seek the CSs within 3 km of the rider's destination. Applying this setting, we further assume the stochastic charging waiting time at a given CS ($WT_{u}^{t}$) is determined by the number of idle EVs within 3 km since we do not have any real charging waiting time data at hand. Specifically, as indicated in~\cite{keskin2021simulation}, the charging time can be drawn from any type of distribution with a known mean and standard deviation. We assume 
the charging time follows a uniform distribution $\mathcal{U} (a, b)$, where the parameters are estimated by the following equations.

\begin{equation}
\begin{aligned}
a = \frac{0.8 - SoC_{max}}{CR}
\end{aligned} 
\end{equation}

\begin{equation}
\begin{aligned}
b = \frac{0.8 - SoC_{min}}{CR}
\end{aligned}
\end{equation}

where $CR$ denotes a constant charging rate, $SoC_{max}$ and $SoC_{min}$ represent the EVs with the highest SoC and lowest SoC, respectively. Further, let $X_1, X_2, \cdots, X_m$ denote the random variables of charging time, the random variable of charging waiting time can be formulated as $WT = X_1 + X_2 + \cdots + X_m$, where $m = \lfloor \frac{|\mathcal{\Bar{D}}|}{n} \rfloor$. $|\mathcal{\Bar{D}}|$ denotes the number of EVs waiting for charging, and $n$ denotes the number of chargers at the given charging station, which is available in the charging station data set. According to the central limit theorem (CLT), the charging waiting time eventually follows a Gaussian distribution given by $\mathcal{N} (\frac{m (a+b)}{2}, \frac{{m}^{2} (b-a)^{2}}{12})$. We acknowledge that the assumption may not hold in the real-world scenario; however, this input parameter does not change the characteristics of the EV batched matching model (the second stochastic optimization in the integrated framework). Hence, it is sufficient to reflect our contribution. Once the real charging waiting time information is available, it can be directly applied to our proposed framework.

Finally, at the beginning of the system operates, idle EVs are randomly generated around the regions as depicted in Fig.~\ref{fig:pick-up}, the initial SoC of EVs is assumed to follow a uniform distribution range in [0,2, 0.8], and three types of EV with different battery consumption rates are assumed to serve rides as described in Table~\ref{table:model-parameter}.  In addition, the fleet size (number of idle EVs) $\mathcal{C}^{t}$ during the batching window $t$ is estimated by trip records in this work. Specifically, we adopt the number of EVs that complete the rider trips at the current batching window as the fleet size.

\subsection{Evaluation Metrics}
We propose three evaluation metrics for the performance comparison from distinct role perspectives: Matching Rate (MR), Rider's Average Waiting Time (RAWT), and Average Charging Waiting Time (ACWT). Firstly, from the standpoint of the EV ride-hailing platform, the MR gauges the system-wide efficiency, which is a key criterion in ride-hailing systems. Secondly, the RAWT serves as a pivotal measurement to quantify the service quality and riders' satisfaction. It quantifies the duration that riders spend waiting for service, with a lower RAWT being indicative of enhanced rider preference for the EV ride-hailing service. Lastly, the ACWT measures the EV drivers' satisfaction, emphasizing the importance of minimal charging waiting times to ensure a consistent supply within the system.

\begin{itemize}
    \item \textbf{Matching Rate (MR)}. The MR of the given batching window is defined as the ratio of the number of matched riders to the number of rider requests, which is shown as follows:

    \begin{equation}
    \begin{aligned}
    MR = \frac{\sum_{j \in \mathcal{D}} \sum_{k \in \mathcal{R}} y_{j, k}}{\left| \mathcal{R} \right|}
    \end{aligned}
    \end{equation}

    \item \textbf{Rider's Average Waiting Time (RAWT)}. The RAWT of the given batching window is defined as the ratio of the sum of rider waiting time to the number of matched riders, which is defined as follows: 

    \begin{equation}
    \begin{aligned}
    RAWT = \frac{\sum_{j \in \mathcal{D}} \sum_{k \in \mathcal{R}} W^{2}_{j, k} \cdot y_{j,k}}{\sum_{j \in \mathcal{D}} \sum_{k \in \mathcal{R}} y_{j, k}}
    \end{aligned}
    \end{equation}
    
    where $W^{2}_{j, k}$ is the rider's waiting time that can be computed in Eq.~(\ref{constraint:c11}).

    \item \textbf{Average Charging Waiting Time (ACWT)}. The ACWT of the given batching window is defined as the ratio of the total EVs' waiting time at CSs to the number of matched EVs, which is defined as follows:

    \begin{equation}
    \begin{aligned}
    ACWT = \frac{\sum_{j \in \mathcal{D}} WT_{j}}{\sum_{j \in \mathcal{D}} \sum_{k \in \mathcal{R}} y_{j, k}}
    \end{aligned}
    \end{equation}
    
    where $WT_{j}$~\footnote{$WT_{u}^{t}$ in Eq.~(\ref{SO-matching-obj:term-2}) denote the random variable of the waiting time in CS $u$ at time $t$ since the exact charging waiting time is unknown during the matching optimization stage, while $WT_{j}$ is a scalar value, which is an outcome of $WT_{u}^{t}$. We adopt the informal notations here for simplicity.} denotes the EV driver $j$'s charging waiting time at CS.
\end{itemize}

\subsection{Performance Evaluation}
This research primarily focuses on evaluating the benefits derived from integrating proactive guidance and CS selection uncertainty within batched matching solutions. To this end, the analysis is structured around two metric dimensions. The first dimension pertains to the guidance strategy, encompassing scenarios of non-guidance, guidance predicated on deterministic demand forecasting, and guidance based on probabilistic demand forecasting to account for uncertainty. The second dimension revolves around batched matching strategies, specifically considering rider waiting time, CS selection, and their combination. In alignment with these dimensions, four benchmark models are selected for comparative analysis. The first model is Batched Matching with Charging Station Selection under No Guidance (BMCSS-NG), which is described in Eq.(~\ref{model-2}). The second model is Batched Matching with Charging Station Selection under Deterministic Guidance (BMCSS-DG), which neglects the forecasting rider demand uncertainty in Eq.(~\ref{model-1}). The third model is Batched Matching with Rider Waiting Time minimization under Stochastic Guidance (BMRWT-SG), which minimizes the rider's average waiting time in the batching matching model. The fourth model is Batched Matching with Charging Waiting Time minimization under Stochastic Guidance (BMCWT-SG), which minimizes the charging waiting time in the batched matching model. Finally, our proposed model is Batched Matching with Charging Station Selection under Stochastic Guidance (BMCSS-SG) as described in the previous section. The selected benchmark models are summarized in Table~\ref{table:benchmarks}. In addition, each group of experiment results are categorized into weekday and weekend scenarios. To ensure a robust long-term perspective, the presented numeric results represent averages derived from six weekdays and six weekends within January 2016. Further, to explain the comparative results more clearly, we display two granularity levels of comparative results. One is the results across daily time slots, which are illustrated in line plots; the other is daily mean results averaged by daily time slots, which are shown in tables. 

\begin{table*}[htbp]
\centering
\small % This sets the font size to small,  \tiny, \scriptsize, \footnotesize,\small, \normalsize, \large, \Large (the 'L' is in uppercase), \LARGE (all uppercase), \huge, \Huge (the 'H' is in uppercase)
\caption{The selected benchmark models for performance comparison, the proposed approach is at bottom right marked in bold}
\renewcommand{\arraystretch}{1.3}  %
\setlength{\tabcolsep}{3.5mm}{
\begin{tabular}{|c|c|c|c|} 
\hline
\diagbox[]{Guidance}{Matching} & Rider Waiting Time & CS Selection & Rider Waiting Time + CS Selection     \\ \hline
No   &  \diagbox[dir=NE]{}{}   & \diagbox[dir=NE]{}{}  & BMCSS-NG \\ \hline
Deterministic & \diagbox[dir=NE]{}{}  &  \diagbox[dir=NE]{}{} & BMCSS-DG \\ \hline
Stochastic & BMRWT-SG & BMCWT-SG & \textbf{BMCSS-SG} \\ \hline
\end{tabular}}
\label{table:benchmarks}
\end{table*}

\subsubsection{Comparative Results of Matching Rate}

First, we compare BMCSS-NG and BMCSS-DG to BMCSS-SG in terms of matching rates. The reason we do not select BMRWT-SG and BMCWT-SG is that the matching rate is influenced by the guidance strategy, which is aligned along the vertical axis of the analytical model, rather than by the matching objective, which is situated along the horizontal axis. Take the stochastic guidance strategy, for example; the set of guided EVs is identical for BMRWT-SG, BMCWT-SG, and BMCSS-SG because all of these models adopt the stochastic guidance strategy. Apparently, from the perspective of optimization, the three models are identical in terms of the feasible regions (the sets of the feasible solutions are identical). Consequently, it can be inferred that the matching rate, as calculated by the models arrayed along the horizontal dimension of Table~\ref{table:benchmarks}, is inherently consistent across these models. This consistency is underpinned by their shared feasible solution sets and uniform application of the stochastic guidance strategy despite the variability in their specific objective functions.

\begin{figure}[!htb]
\minipage{0.5\textwidth}
  \includegraphics[width=\linewidth]{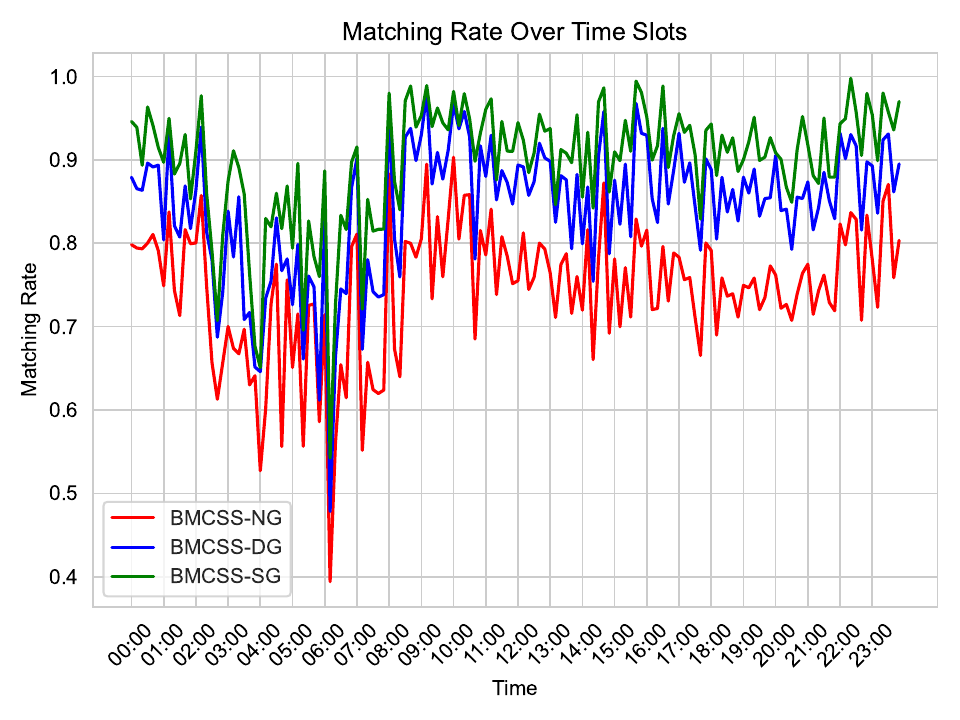}
  \caption{Matching Rate Comparison (weekdays)}\label{fig:MR-weekday}
\endminipage\hfill
\minipage{0.5\textwidth}
  \includegraphics[width=\linewidth]{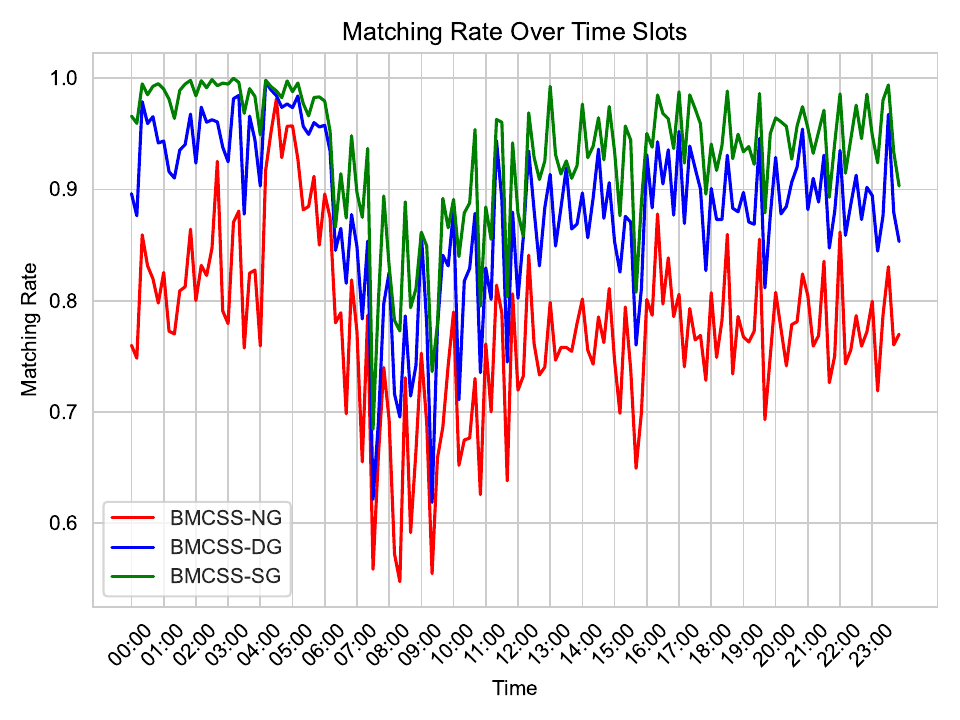}
  \caption{Matching Rate Comparison (weekends)}\label{fig:MR-weekend}
\endminipage\hfill
\end{figure}

The comparative results of matching rates between BMCSS-SG and benchmark models on the weekdays and weekends are visualized in Fig.~\ref{fig:MR-weekday} and \ref{fig:MR-weekend}, respectively. It can be observed that compared to BMCSS-NG, the matching rate from BMCSS-SG is much higher, and the trend is consistent over the daily time slots in both weekday and weekend scenarios. This is because the available idle EVs are guided to the designated regions around the POI, which implies that the idle EVs are close to the potential riders' origins. Hence, more EVs are available to meet the rider requests, which leads to a higher matching rate. In addition, BMCSS-SG also demonstrates superior performance in matching rates compared to BMCSS-DG. This is because, compared to the point prediction, the probabilistic forecasting results can provide informative rider demand information, which enhances the stochastic optimization model to obtain better solutions. This illustrates the intrinsic benefit of incorporating rider demand uncertainty within the optimization process, as evidenced by the improved matching rates.

Further, the summarized comparison of matching rates under varied guidance strategies is presented in Table~\ref{table:MR}. It can be observed that BMCSS-SG are 6.1\% and 20.1\% on average, up to 8.9\% and 22.9\% higher than BMCSS-DG and BMCSS-NG in terms of the matching rate during the weekdays. Analogous trends are observed during weekend periods, where the average increments in matching rates for BMCSS-SG are 5.4\% and 19.9\% relative to BMCSS-DG and BMCSS-NG, with the maximum increases being 7.8\% and 22.1\%, respectively. These findings underscore the significant advantages that an EV ride-hailing platform can accrue through the adoption of a stochastic guidance strategy, as evidenced by the marked improvements in matching rates. This enhancement is indicative of the efficacy of incorporating probabilistic elements into the guidance strategy, thereby optimizing the allocation of EV resources in response to fluctuating demand patterns.

\begin{table}[]
\centering
\small % This sets the font size to small,  \tiny, \scriptsize, \footnotesize,\small, \normalsize, \large, \Large (the 'L' is in uppercase), \LARGE (all uppercase), \huge, \Huge (the 'H' is in uppercase)
\caption{The statistics of the daily matching rate (MR), including mean, maximum, minimum, and standard deviation averaged by six-day results under different guidance strategies, the highest values are marked in bold}
\renewcommand{\arraystretch}{1.3}  %
\setlength{\tabcolsep}{5.5mm}{
\begin{tabular}{c|c|cccc}
\toprule[1pt]
\multicolumn{2}{c|}{\diagbox[]{Models}{MR}} & Mean & Max & Min & SD \\ 
\hline
Weekdays & BMCSS-NG           &  74.52\%    &  74.79\%   &  74.09\%   &  0.38\%   \\
        & BMCSS-DG &   84.35\%   & 85.33\%    &   83.62\%   &   0.74\%  \\
        & BMCSS-SG   &  \textbf{89.52\%}    &  \textbf{91.06\%}   &  \textbf{88.97\%}   &   1.03\%  \\
\hline
Weekends & BMCSS-NG           &   77.65\%   &  78.39\%   &  76.97\%   &  0.71\%   \\
        & BMCSS-DG &   88.33\%    &  88.87\%   &  87.16\%   &   0.79\%  \\
        & BMCSS-SG   &  \textbf{93.10\%}    &   \textbf{93.97\%}  &   \textbf{91.81\%}  &  0.92\%   \\
\bottomrule
\end{tabular}}
\label{table:MR}
\end{table}

\subsubsection{Comparative Results of Rider's Average Waiting Time}

\begin{figure}[!htb]
\minipage{0.5\textwidth}
  \includegraphics[width=\linewidth]{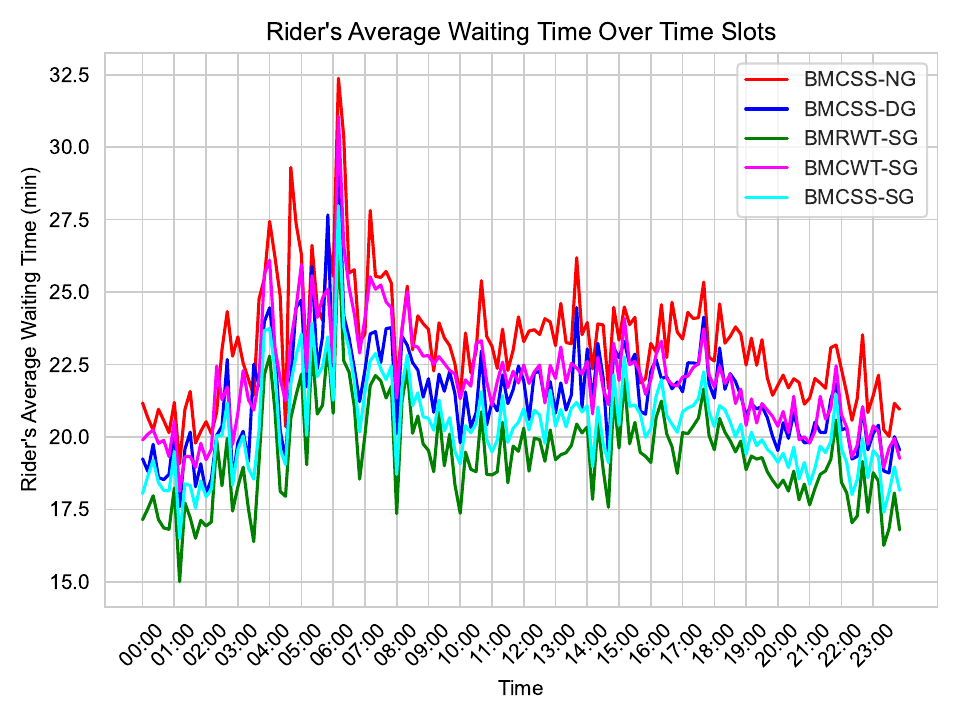}
  \caption{Rider Average Waiting Time Comparison (weekdays)}\label{fig:RWT(weekday)}
\endminipage\hfill
\minipage{0.5\textwidth}
  \includegraphics[width=\linewidth]{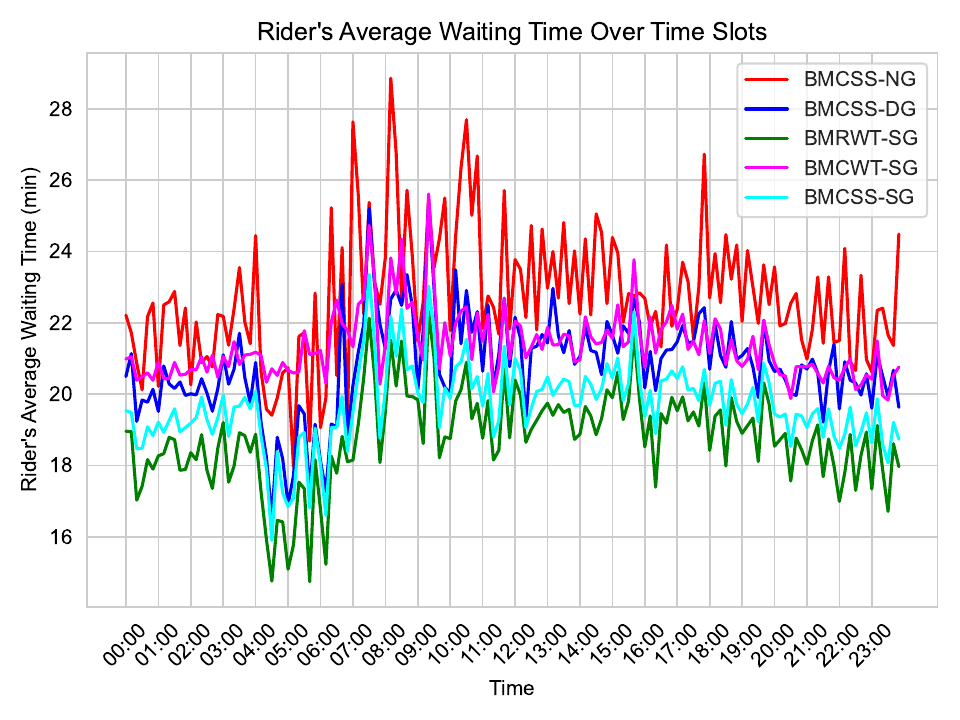}
  \caption{Rider Average Waiting Time Comparison (weekends)}\label{fig:RWT(weekend)}
\endminipage\hfill
\end{figure}

Second, we compare BMCSS-NG, BMCSS-DG, BMRWT-SG, and BMCWT-SG to BMCSS-SG in terms of rider's average waiting time. The comparison results across daily time slots in the weekdays and weekends are illustrated in Fig.~\ref{fig:RWT(weekday)} and Fig.~\ref{fig:RWT(weekend)}, respectively. Similar to the results in the matching rate comparison, it is observed that our model BMCSS-SG outperforms BMCSS-NG and BMCSS-DG. As indicated in Table~\ref{table:RAWT}, in comparison to BMCSS-NG and BMCSS-DG, the rider's average waiting time experienced by riders is observed to decrease by 12.7\% and 4.2\% during weekdays, with maximum reductions noted at 14.3\% and 5.8\%. Similar trends are observed during weekends, with average reductions in waiting time being 13.1\% and 5.2\%, and the peak reductions reaching 14.9\% and 6.6\%. This reduction in waiting time is attributed to two primary factors. Firstly, the strategic deployment of idle EVs to regions proximal to rider trip origins significantly narrows the physical distance between EVs and riders, resulting in decreased waiting times. Secondly, as discussed in the preceding matching rate results, an augmented matching rate inherently contributes to a reduction in the average waiting time for riders. This is based on the assumption that riders not matched within the existing batching window are required to wait until the consequent time slot for potential matching. While this assumption may appear to be stringent, it aligns with the core characteristics of the problem framework. In a more pragmatic scenario, riders might choose to cancel their initial request via the platform and subsequently resubmit a new request after a period. This real-world scenario is effectively analogous to the assumption utilized in this study, thereby reinforcing the validity of the chosen evaluation metric in reflecting the operational efficiency of the EV ride-hailing platform.

Furthermore, as expected, the proposed model BMCSS-SG is between the BMRWT-SG and BMCWT-SG in terms of the rider's average waiting time. This observation is logically coherent, considering that BMRWT-SG is explicitly oriented towards minimizing rider waiting time in its objective function, whereas BMCWT-SG is focused on the minimization of charging waiting time. While BMCSS-SG seeks the balance between the two types of waiting time minimization. Notably, the application of BMCSS-SG results in the rider's average waiting time being marginally higher, specifically by 5.3\% and 5.1\% in weekdays and weekends, respectively, compared to BMRWT-SG. Although BMRWT-SG is optimized for the lowest rider waiting time, it does not account for charging waiting time, which is a crucial factor impacting the supply dynamics within the EV ride-hailing platform. In this context, BMCSS-SG emerges as a balanced, trade-off solution, effectively harmonizing rider waiting time with charging waiting time, thus catering to the interests of both EV drivers and riders. Consequently, BMCSS-SG presents itself as a more pragmatic model for the EV ride-hailing market. Moreover, the BMCWT-SG is slightly inferior to BMCSS-DG. The outcome makes sense because BMCSS-DG also incorporates rider waiting time as part of the objective function. It is worth noting that despite the exclusion of rider waiting time in its objective framework, BMCWT-SG still surpasses BMCSS-NG regarding the average waiting time for riders. This implies that the strategic guidance operations for EVs can mitigate the impact of not including rider waiting time in the objective function.

\begin{table}[]
\centering
\small % This sets the font size to small,  \tiny, \scriptsize, \footnotesize,\small, \normalsize, \large, \Large (the 'L' is in uppercase), \LARGE (all uppercase), \huge, \Huge (the 'H' is in uppercase)
\caption{The statistics of the daily rider's average waiting time (RAWT), including mean, maximum, minimum, and standard deviation averaged by six-day results under different benchmark models (unit in minutes), the lowest values are marked in bold}
\renewcommand{\arraystretch}{1.3}  %
\setlength{\tabcolsep}{5.5mm}{
\begin{tabular}{c|c|cccc}
\toprule[1pt]
\multicolumn{2}{c|}{\diagbox[]{Models}{RAWT}} & Mean & Max & Min & SD \\ 
\hline
Weekdays & BMCSS-NG           &  23.26    &   23.53  &  22.94   &  0.31   \\
        & BMCSS-DG &  21.32    &  21.40  &   21.11  &  0.14   \\
        & BMRWT-SG   &  \textbf{19.41}    &  \textbf{19.66}   &  \textbf{19.20}   &  0.20   \\
        & BMCWT-SG &  22.01    &  22.32   &   21.79   &  0.22   \\
        & BMCSS-SG   &  20.43    & 20.69    &  20.17   &   0.21  \\
\hline
Weekends & BMCSS-NG           &  22.72    &  22.95   &  22.45   &  0.10   \\
        & BMCSS-DG &   20.82   &  20.91  &  20.77   &  0.06   \\
        & BMRWT-SG   &  \textbf{18.78}    &  \textbf{18.95}   &  \textbf{18.53}   &   0.19  \\
        & BMCWT-SG &  21.35    &  21.52   &   21.12   &   0.17  \\
        & BMCSS-SG   &  19.74    &  19.87   &  19.53   &   0.15  \\
\bottomrule
\end{tabular}}
\label{table:RAWT}
\end{table}

\subsubsection{Comparative Results of Average Charging Waiting Time}

Finally, we compare BMRWT-SG and BMCWT-SG to BMCSS-SG in terms of charging waiting time. The reason we select the benchmark models along the horizontal direction in Table~\ref{table:benchmarks} stems from an interest in discerning the extent of benefits that an EV ride-hailing platform may accrue from integrating CS selection decisions, specifically assessing the potential reduction in charging waiting time for EVs with low SoC. In this sense, comparing the charging waiting time under different set of guided EVs is meaningless, therefore, we compare BMCSS-SG to BMRWT-SG and BMCWT-SG, because the set of guided EVs are identical under the same guidance strategy. To ensure a relevant comparison, the analysis focuses on two distinct categories of use cases: one is the charging waiting time of EVs with low SoC, defined as being less than 30, and the other is the charging waiting time of EVs with medium SoC, quantified as being between 30\% and 60\%. It is critical to reiterate that the precise modeling of exact charging waiting times is beyond the scope of this study. In the following experimental setup, the charging waiting time at each CS is approximated based on the count of idle EVs (those having completed their rider drop-offs), a metric that may not entirely align with real-world scenarios. Nonetheless, it is important to note that such an approximation does not detract from the validity of the comparative analysis between the proposed model and the benchmark models. This methodology provides valuable insights into the efficacy of the integrated CS selection strategy, particularly in terms of its impact on charging waiting times under varying SoC conditions.

\begin{figure}[!htb]
\minipage{0.5\textwidth}
  \includegraphics[width=\linewidth]{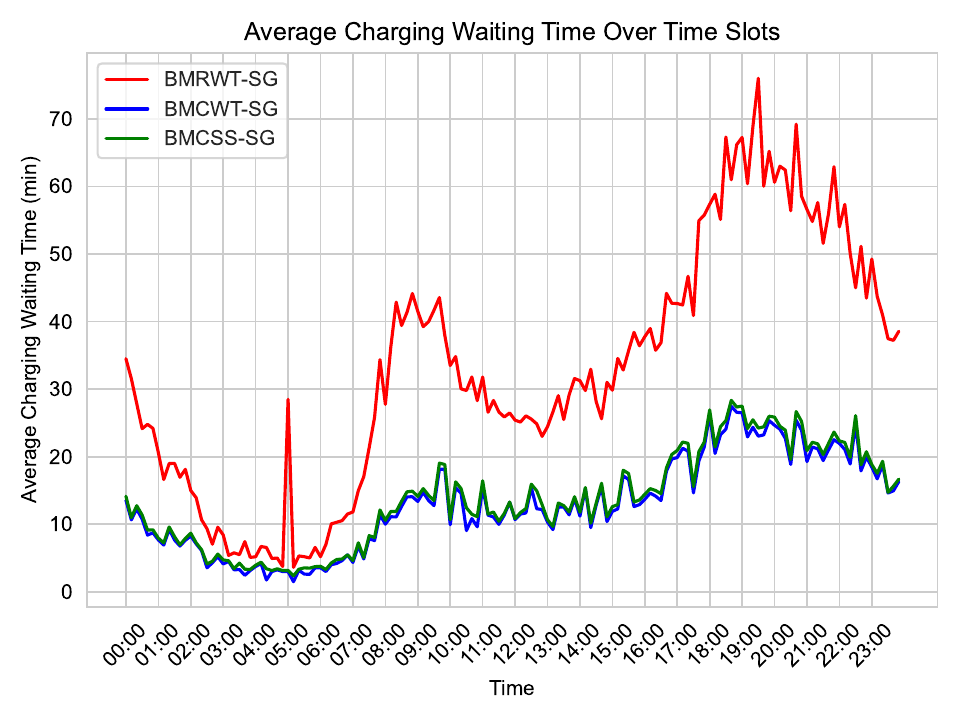}
  \caption{Average Charging Waiting Time of Low SoC (weekdays)}\label{fig:CWT-30 (weekday)}
\endminipage\hfill
\minipage{0.5\textwidth}
  \includegraphics[width=\linewidth]{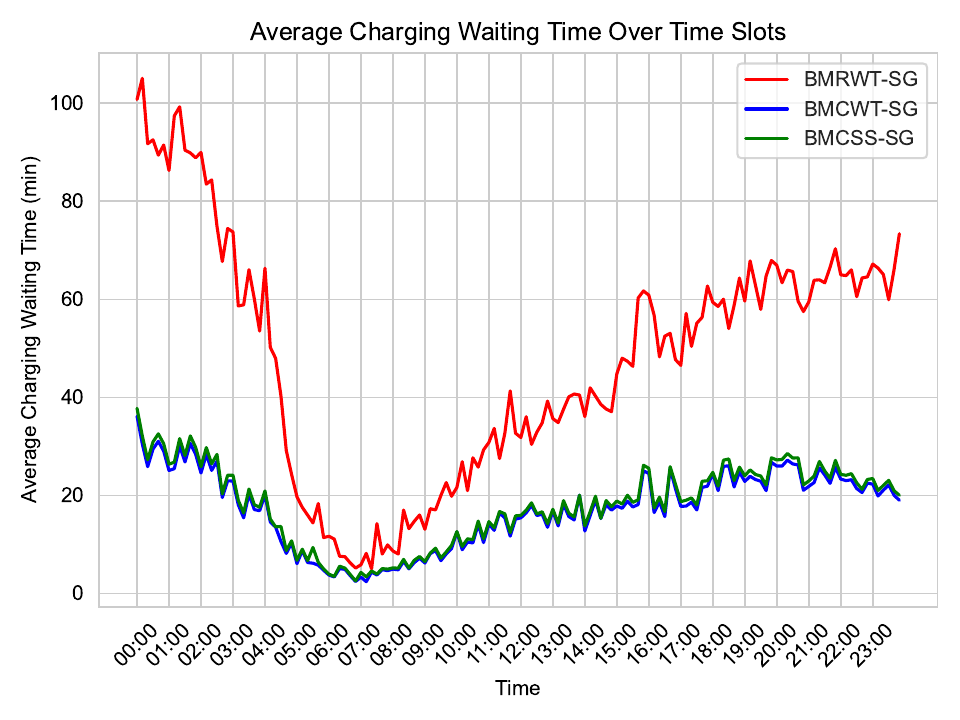}
  \caption{Average Charging Waiting Time of Low SoC (weekends)}\label{fig:CWT-30 (weekend)}
\endminipage\hfill
\end{figure}

The average charging waiting time of EVs with low SoC (less than 30\%) across daily time slots on weekdays and weekends are illustrated in Fig.~\ref{fig:CWT-30 (weekday)} and~\ref{fig:CWT-30 (weekend)}. Observations reveal that both BMCWT-SG and BMCSS-SG demonstrate a consistently superior performance over the BMRWT-SG in terms of average charging waiting time across all the time slots. Specifically, the average charging waiting time for BMRWT-SG is found to be 151.38\% and 175.68\% greater than that of BMCSS-SG during weekdays and weekend scenarios, respectively. This is because BMRWT-SG exclusively prioritizes minimizing rider waiting time, without integrating CS selection decisions into its objective function. As a result, EVs with low SoC might be allocated to riders whose destinations are not in proximity to CSs with shorter charging wait times, leading to prolonged charging waiting times. Additionally, it is noteworthy that the charging waiting time for BMCSS-SG is marginally higher than BMCWT-SG, by 4.4\% and 4.7\% on weekdays and weekends, respectively. This is attributed to the fact that BMCSS-SG endeavors to simultaneously minimize both rider waiting time and charging waiting time, which inadvertently results in a somewhat increased charging waiting time compared to BMCWT-SG. However, as previously discussed, BMCWT-SG incurs a higher average rider waiting time in comparison to BMCSS-SG, which could potentially compromise the rider experience. From this perspective, BMCSS-SG emerges as a balanced, trade-off solution, effectively harmonizing charging waiting time with rider waiting time. This trade-off indicates that the model offers mutual benefits to both EV drivers and riders, endorsing its utility as a comprehensive solution in the context of EV ride-hailing services.

\begin{figure}[!htb]
\minipage{0.5\textwidth}
  \includegraphics[width=\linewidth]{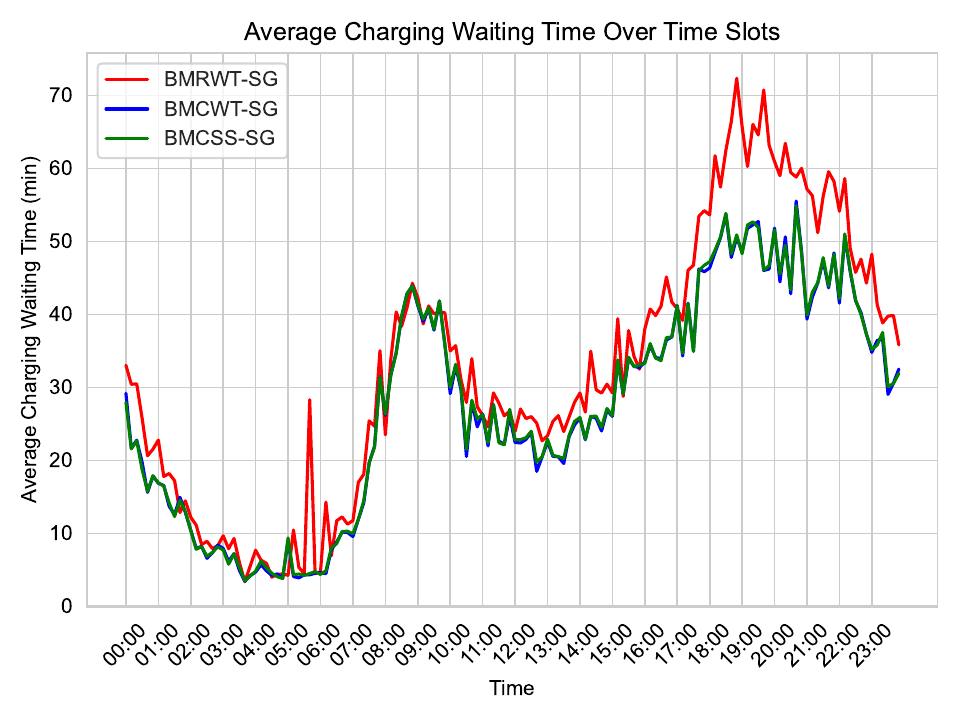}
  \caption{Average Charging Waiting Time of Medium SoC (weekdays)}\label{fig:CWT-60 (weekday)}
\endminipage\hfill
\minipage{0.5\textwidth}
  \includegraphics[width=\linewidth]{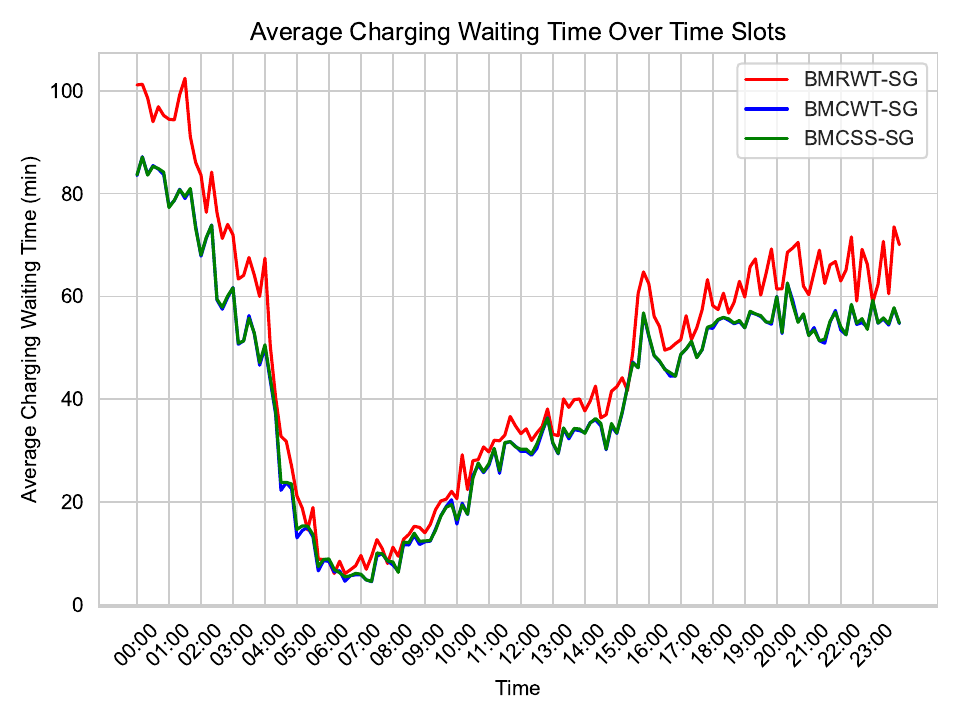}
  \caption{Average Charging Waiting Time of Medium SoC (weekends)}\label{fig:CWT-60 (weekend)}
\endminipage\hfill
\end{figure}

The average charging waiting time of EVs with medium SoC are plotted in Fig.~\ref{fig:CWT-60 (weekday)} and~\ref{fig:CWT-60 (weekend)}. Similar to the comparative results in the low SoC scenario, both BMCWT-SG and BMCSS-SG exhibit superior performance in reducing charging waiting times when contrasted with BMRWT-SG. However, the difference in performance is less significant in the medium SoC scenario compared to the low SoC scenario. More precisely, the discrepancy in charging waiting times, wherein BMRWT-SG exceeds BMCWT-SG, is quantified at 18.69\% and 16.88\% for weekdays and weekends, respectively. This variation stands in stark contrast to the more substantial differences observed in the low SoC scenario (151.38\% and 175.68\%, as previously noted). Remarkably, there exist specific time slots (for instance, 8:00 on weekdays) where BMRWT-SG demonstrates a shorter charging waiting time compared to both BMCWT-SG and BMCSS-SG. This phenomenon can be attributed to the allocation of EVs with higher SoC to rider destinations associated with longer charging waiting times. Clearly, the implementation of BMCWT-SG and BMCSS-SG ensures that EVs with lower SoC are preferentially directed to the CSs characterized by shorter charging waiting times. Nonetheless, it is imperative to acknowledge that BMCWT-SG achieves this objective at the expense of rider satisfaction, evident through extended rider's waiting times. Such a compromise may potentially be deemed untenable from the perspective of the EV ride-hailing platform, which underscores the importance of balancing both charging waiting times and rider satisfaction within these operational models.

In summary, the comparative results suggest that our decision-making framework BMCSS-SG tackles the EV batched matching with the CS selection problem in a more reasonable manner. In fact, this trade-off solution can be controlled by calibrating the parameters (i.e., $\theta_{1}$ and $\theta_{2}$) in the model's objective function. For instance, in scenarios typified by high traffic volumes, such as during rush hours, the system operator might prioritize $\theta_{2}$ over $\theta_{1}$ to ensure minimal rider waiting times. Conversely,  in periods characterized by lower demand, such as during off-peak hours, $\theta_{1}$ may be accentuated over $\theta_{2}$ to facilitate shorter charging waiting times. This parameter adjustment capability imbues the BMCSS-SG model with a heightened degree of adaptability, thereby rendering it a particularly versatile and dynamic solution for EV ride-hailing services. The model's inherent flexibility to adapt to varying operational conditions positions it as a robust framework capable of optimizing service efficiency under diverse circumstances.

\begin{table}[]
\centering
\small % This sets the font size to small,  \tiny, \scriptsize, \footnotesize,\small, \normalsize, \large, \Large (the 'L' is in uppercase), \LARGE (all uppercase), \huge, \Huge (the 'H' is in uppercase)
\caption{The statistics of the daily average charging waiting time (ACWT), including mean, maximum, minimum, and standard deviation averaged by six-day results, under different benchmark models (unit in minutes), the lowest values are marked in bold}
\renewcommand{\arraystretch}{1.3}
\setlength{\tabcolsep}{5.5mm}{
\begin{tabular}{c|c|c|cccc}
\toprule[1pt]
\multicolumn{3}{c|}{\diagbox[]{Models}{ACWT}} & Mean & Max & Min & SD \\ 
\hline
\multirow{6}{*}{SoC $\leq$ 30\%} & Weekdays & BMRWT-SG   &  32.78   &  36.81   & 30.87  &   2.76 \\
        & & BMCWT-SG &  \textbf{13.04}    &  \textbf{14.54}  & \textbf{12.37}   &  1.01  \\
        & & BMCSS-SG   &  13.61    & 15.12  &  12.94   &  1.07  \\
\cline{2-7}
& Weekends & BMRWT-SG   & 47.72     & 55.48   &  37.15  &  7.91   \\
        & & BMCWT-SG &  \textbf{17.31}   & \textbf{19.27}    &  \textbf{14.06}    &  2.42  \\
        & & BMCSS-SG   &  18.07   & 20.24    &   14.65  &  2.56  \\

\hline
\multirow{6}{*}{30\% $\leq$ SoC $\leq$ 60\%} & Weekdays & BMRWT-SG   &  32.52   &  35.67   &  31.02 &  2.17  \\
        & & BMCWT-SG &  \textbf{27.40}    &  \textbf{30.86}  &   \textbf{26.02}  &  2.31  \\
        & & BMCSS-SG   & 27.53     & 30.98  &  26.18   &  2.31  \\
\cline{2-7}
& Weekends & BMRWT-SG   & 48.19     &  55.66  & 37.09   &  8.12   \\
        & & BMCWT-SG &  \textbf{41.23}   &  \textbf{48.36}   &  \textbf{31.22}    &   7.41 \\
        & & BMCSS-SG   &  41.41   &  48.46   &   31.38  &  7.39  \\
\bottomrule
\end{tabular}}
\label{table:ACWT}
\end{table}

\section{Conclusions and Future Work}\label{section:conclusion}
In this study, we presented a holistic stochastic optimization framework designed to tackle the multi-faceted challenges inherent in EV ride-hailing services. Our approach, addressing uncertainties in rider demand and CS selection, integrates EV matching, proactive guidance, and CS selection into a unified data-driven optimization framework. Through extensive numeric experiments, we have demonstrated the efficacy of our framework in enhancing the efficiency and adaptability of EV ride-hailing operations. We believe that the EV ride-hailing platform can be beneficial from this flexibly integrated solution. It not only offers immediate practical benefits by improving the operational efficiency of EV ride-hailing services but also provides a foundation for further research and innovation in sustainable urban mobility. 

We identify two research directions for our future work. Firstly, charging station selection is determined by two factors (EV travel distance and charging waiting time) during the batched matching optimization in this work. Nonetheless, an EV driver might opt to continue providing services for subsequent rides rather than engaging in recharging activities, contingent upon sufficient battery levels, high charging demand, or escalated rider demand. In this sense, the driver behavior can be modeled by discrete choice and/or multinomial logit model to make a more reasonable charging decision. Secondly, it is acknowledged that the EV ride-hailing platform may incur profit losses when the available supply fails to meet rider demand. The concept of surge pricing in charging stations emerges as a critical element in reducing charging waiting times to alleviate the problem. Despite its significance, the extent of financial expenditure that the EV ride-hailing platform is prepared to incur under the surge pricing model remains unresolved.


% \bibliographystyle{unsrt}
% \bibliographystyle{plainnat}
% \bibliographystyle{elsarticle-num}
% \bibliographystyle{elsarticle-harv}
% \bibliographystyle{elsarticle-num-names}

\bibliographystyle{elsarticle-num} 
\bibliography{cas-refs}





% \end{thebibliography}
\end{document}